\documentclass[10pt]{article}
\usepackage{amsmath,amssymb,amsfonts}

\newtheorem{theorem}{Theorem}

\newcommand{\ZZ}{\mathbb{Z}}
\newcommand{\Sym}{\rm Sym}

\begin{document}

\title{A note on 2-subset-regular self-complementary 3-uniform hypergraphs}
\author{Martin Knor\footnote{Department of Mathematics, 
Faculty of Civil Engineering, Slovak University of Technology, 
Radlinsk\'eho 11, 813 68 Bratislava, Slovakia, {\tt knor@math.sk}.}\hskip2mm
and
Primo\v{z} Poto\v cnik\footnote{Faculty of Mathematics and Physics,
University of Ljubljana, Slovenia, {\sc and} Institute of Mathematics, Physics, 
and Mechanics, Jadranska 19, Sl-1000 Ljubljana, Slovenia, 
{\tt primoz.potocnik@fmf.uni-lj.si}.}}

\date{}\maketitle

\begin{abstract}
We show that a 2-subset-regular self-complementary 3-uniform hypergraph
with $n$ vertices exists if and only if $n\ge 6$ and $n$ is congruent to $2$
modulo $4$.
\end{abstract}

\section{Introduction}

A {\em $k$-uniform hypergraph} of order $n$ is an ordered pair $\Gamma = (V,E)$,
where $V = V(\Gamma)$ is an arbitrary set of size $n$, and $E=E(\Gamma)$ is a subset of
$V^{(k)} = \{ e\subseteq V : |e| = k\}$.
Note that the notion of a $2$-uniform hypergraph coincides with the usual
notion of a simple graph. We shall call a $k$-uniform hypergraph simply a {\em $k$-hypergraph}.

A $k$-hypergraph $\Gamma$ is {\em self-complementary} if it is
isomorphic to its complement $\Gamma^C$, defined by $V(\Gamma^C) = V(\Gamma)$
and $E(\Gamma^C) = V(\Gamma)^{(k)}\setminus E(\Gamma)$. Equivalently,
$\Gamma = (V,E)$ is self-complementary whenever there exists a permutation $\tau\in \Sym(V)$, called the
{\em antimorphism of} $\Gamma$, such that for all $e\in V^{(k)}$ the equivalence 
$e\in E \Leftrightarrow e^\tau \not \in E(\Gamma)$ holds. 
Antimorphisms of uniform hypergraphs were characterized 
in terms of their cyclic decompositions by Wojda in \cite{wojda}.

A $k$-hypergraph $\Gamma$ is {\em $t$-subset-regular} if there exists an integer $\lambda$,
also called the {\em $t$-valence} of $\Gamma$, such that each $t$-element subset of $V(\Gamma)$ 
is a subset of exactly $\lambda$ elements of $E(\Gamma)$. Clearly $t$-subset-regular 
$k$-hypergraphs generalize the notion of regular graphs,
and can also be viewed as a bridge between graph theory and design theory. Namely, 
a $t$-subset-regular $k$-hypergraph of $t$-valence $\lambda$ and order $n$ is 
simply a $t$-$(n,k,\lambda)$ design.
Moreover, such a $k$-hypergraph $\Gamma$ is self-complementary if and only if the pair $\{\Gamma,\Gamma^C\}$
forms a large set of $t$-designs LS$[2](t,n,k)$ with the additional property that the two designs
constituting the large set are isomorphic (see \cite{hb} for the definition of a large set of designs).

Here a question of existence of a self-complementary $t$-subset-regular hypergraph with prescribed parameters
$n$, $k$, and $t$ arises naturally. An easy counting argument shows that whenever a self-complementary
$t$-subset-regular $k$-hypergraphs on $n$ vertices exists, then
${n-i\choose k-i}$ is even for all $i=0, \ldots, t$. 

It can be seen that the above divisibility conditions can in fact be expressed in terms of certain congruence
conditions on $n$ modulo an appropriate power of $2$ (see \cite{mod}). 
For example if $k = 2^\ell$ or $k = 2^\ell + 1$ for some positive integer $\ell$, then
$n$ is congruent to one of $t, \ldots, k-1$ modulo $2^{\ell+1}$. In particular, if $k=2$ and $t=1$,
then $n \equiv 1\,\, (\hbox{mod } 4)$;  if $k=3$ and $t=1$, then $n \equiv 1$ or $2\,\,
(\hbox{mod } 4)$; if $k=3$ and $t=2$, then $n \equiv 2\,\, (\hbox{mod } 4)$.

In \cite{PS} the following question, strengthening Hartman's conjecture \cite{hartman}
about existence of large sets of (not necessarily isomorphic) designs, was raised:
\medskip

{\bf Question.} \cite{PS} {\em Is it true that for every triple of integers $t<k<n$ such that
${n-i\choose k-i}$ is even for all $i=0, \ldots, t$, there exists
a self-complementary $t$-subset-regular $k$-hypergraph of order $n$?}
\medskip

It is well known (see \cite{sachs}) that a regular self-complementary graph on $n$ vertices 
exists if and only if $n$ is congruent to $1$ modulo $4$, showing that the 
answer to the above question is affirmative for $k=2$ and $t=1$.
Recently, the answer was proved to be affirmative also for the case $k=3$ and $t=1$ (see \cite{PS}).
The aim of this note is to show that the answer to the question above is affirmative also for the
remaining case of $3$-hypergraphs, namely for the case $k=3$, $t=2$. More precisely, in Section~\ref{sec:con}
we present a construction which proves the following:

\begin{theorem}
If $n\ge 6$ and $n$ is congruent to $2$ modulo $4$, then there exists
a 2-subset-regular self-complementary $3$-hypergraph on $n$ vertices.
\end{theorem}

\section{Construction}
\label{sec:con}

Let $n=4k+2$ for some integer $k$. For 
$i=0,1$, let $V_i=\{0_i,1_i,\dots,(2k)_i\}$ be a copy of the
the ring $\ZZ_{2k+1}$. Define $\Gamma_n$ to be the $3$-hypergraph with
vertex set $V=V_0\cup V_1$ and edge set
$E=E_1\cup E_2\cup E_3$, where
\begin{eqnarray*}
E_1 &=& V_0^{(3)},\\
E_2 &=& \{\{a_0,b_0,c_1\}\> :\> a,b\in \ZZ_{2k+1},\> a\not = b,\> c=\frac{a+b}{2}\},\\
E_3 &=& \{\{a_0,b_1, c_1\}\> :\> a,b,c\in \ZZ_{2k+1},\> a\not=\frac{b+c}{2}\}.
\end{eqnarray*}
Note that $2$ is invertible in $\ZZ_{2k+1}$, hence dividing by $2$ in the definitions of $E_2$ and $E_3$ is
well defined.

First we show that $\Gamma_n$ is 2-subset-regular, i.e. we show that
each pair of vertices is contained in exactly $(n-2)/2 = 2k$ edges.
There are four types of pairs of vertices to consider:
\vskip 2mm

(a)
A pair $a_0,b_0$, where $a,b\in \ZZ_{2k+1}$, $a\ne b$.
This pair is contained in $2k{-}1$ edges of $E_1$ and in a unique edge in $E_2$.
As it is contained in none of the edges of $E_3$, the pair is in total of $2k$
edges.

(b)
A pair $a_1,b_1$, where $a,b\in \ZZ_{2k+1}$, $a\ne b$.
This pair appears only in $(2k+1) - 1 = 2k$ edges of $E_3$.

(c)
A pair $a_0,a_1$, where $a\in \ZZ_{2k+1}$.
This pair also appears only in the edges of $E_3$. In fact, it appears
precisely in the $2k$ edges of the form $\{a_0,a_1,b_1\}$, $b\in \ZZ_{2k+1}\setminus\{a\}$.

(d)
A pair $a_0, c_1$, where $a,c\in \ZZ_{2k+1}$, $a\ne c$.
This pair is contained in the edge $\{a_0, c_1, (2c-a)_0\}$ of $E_2$,
and in the $2k-1$ edges of the form $\{a_0,c_1,b_1\}$ where $b\in \ZZ_{2k+1} \setminus \{2a-c\}$,
of $E_3$. Hence this pair is in exactly $2k$ edges of $\Gamma_n$.
\vskip 2mm

This proves that $\Gamma_n$ is a 2-subset-regular hypergraph.
To prove that it is self-complementary, note that the mapping $\phi:V\to V$ defined by
$\phi(a_i)=a_{i+1}$, for $i=0,1$, with addition in the subscript being modulo 2, is
an antimorphism of $\Gamma_n$.
\hfill{\framebox[2.5mm]{}}
\vskip 5mm

We remark that $\Gamma_n$ is not a vertex-transitive hypergraph if $n>6$
(check the complete sub-hypergraphs of order $n/2$).
However, in $\Gamma_6$ every pair of vertices appears in exactly two edges.
Hence $\Gamma_6$ can be considered as a triangular embedding of a complete 
graph $K_6$ into a surface (see \cite{GT} for a detailed account on graph embeddings).
In fact, $\Gamma_6$ represents a regular triangulation of the projective
plane by $K_6$.
As a consequence, $\Gamma_6$ is vertex-transitive.

\vskip 1pc
{\bf Acknowledgement.}
This paper was partially supported by the Slovak-Slovenian bilateral project, grant no.\ SK-SI-01906.
Martin Knor acknowledges partial support by the Slovak research grants
VEGA 1/0489/08, APVT-20-000704 and APVV-0040-06.

\end{document}